\documentclass[12pt]{amsart}

\usepackage{amsmath,amssymb,enumerate}
\usepackage{relsize}%

\usepackage{epsfig,fancyhdr,color}

\usepackage{amssymb}
\usepackage{amsmath,amsthm}
\usepackage{bbm}
\usepackage{latexsym}
\usepackage{amscd}
\usepackage{psfrag}
\usepackage{graphicx}
\usepackage[all]{xy}

\usepackage{verbatim}
\usepackage{amssymb}
\usepackage{amsmath}
\usepackage{graphicx}
\usepackage{color}
\usepackage{amsthm}
\usepackage{tikz}
\usetikzlibrary{arrows,positioning,decorations.pathmorphing, decorations.markings}

\usepackage{rotating}
\usepackage{adjustbox}

\usepackage{xcolor}

\usepackage[font=small,labelfont=bf]{caption}

\usepackage{mathtools}

\usepackage{tikz-cd}

\usepackage{stackrel}

\usepackage[utf8]{inputenc}
\usepackage{amsmath}
\usepackage{amssymb}

\usepackage{mathabx,epsfig}
\def\acts{\mathrel{\reflectbox{$\righttoleftarrow$}}}

\usepackage[utf8]{inputenc}

\numberwithin{equation}{section}

\definecolor{NoteColor}{rgb}{1,0,0}

\newtheorem*{theorem 1}{\rm\bf Proposition 1}
\newtheorem*{theorem 2}{\rm\bf Proposition 2}

\theoremstyle{definition}

\theoremstyle{remark}

\def\interieur#1{\mathord{\mathop{\kern 0pt #1}\limits^\circ}}

\def\interieur#1{\mathord{\mathop{\kern 0pt #1}\limits^\circ}}

\def\hyperp{{\rm I}\kern-.3ex{\rm H}}

\usepackage{mathabx,epsfig}
\def\acts{\mathrel{\reflectbox{$\righttoleftarrow$}}}

\newtheorem*{theorem*}{Theorem}

\begin{document}

\title[Finite presentation of P(SL(2,${\mathbb Z}$))]{Universal Spin Teichm\"uller Theory, II. \\Finite presentation of P(SL(2,${\mathbb Z}$))}
\
\author{Robert Penner} \address {\hskip -2.5ex Institut des Hautes \'Etudes Scientifiques\\ 35 route des Chartres\\ Le Bois Marie\\ 91440 Bures-sur-Yvette\\ France\\ {\rm and}~Mathematics Department, UCLA\\ Los Angeles, CA 90095\\USA} \email{rpenner{\char'100}ihes.fr}


 \date{\today}


\begin{abstract} 
In previous works, the universal mapping class group was taken to be the group ${\rm PPSL}(2,{\mathbb Z})$ of all piecewise ${\rm PSL}(2,{\mathbb Z})$ homeomorphisms of the unit circle $S^1$ with finitely many breakpoints among the rational points in $S^1$,  and in fact,
the Thompson group $T\approx{\rm PPSL}(2,{\mathbb Z})$.  The new spin mapping class group P(SL(2,${\mathbb Z}$))
is given by all piecewise-constant maps $S^1\to{\rm SL}(2,{\mathbb Z})$ which 
projectivize to an element of ${\rm PPSL}(2,{\mathbb Z})$.
We compute a finite presentation of ${\rm PPSL}(2,{\mathbb Z})$ 
from basic principles of general position as an orbifold fundamental group.
The orbifold deck group of the spin cover 
is explicitly computed here, from which follows also a finite presentation of ${\rm P(SL(2},{\mathbb Z}))$.  This is our
main new achievement.
Certain commutator relations 
in ${\rm P(SL(2},{\mathbb Z}))$ seem to organize according to root lattices,
which would be a novel development.
We naturally wonder what is the automorphism group of ${\rm P(SL(2},{\mathbb Z}))$
and speculate that it is a large sporadic group.  
There is a companion paper to this one which explains
the topological background from first principles, proves that the group studied here using combinatorial group theory
is indeed ${\rm P(SL(2},{\mathbb Z}))$.
\end{abstract}

\maketitle

\centerline {\it For Valentin Po\'enaru on the}

\centerline{\it  happy occasion of his 90th birthday.}



\setcounter{footnote}{0}

\section*{Introduction}

Consider the space of all tesselations $\tau$ of the Poincar\'e disk ${\mathbb D}$, i.e., $\tau$ is a locally-finite collection of geodesics which decompose ${\mathbb D}$ into ideal triangles, and let us choose  a distinguished oriented edge or {\it doe} in $\tau$.  Let 
${\mathcal Tess}$ be the space of PSL(2,${\mathbb R}$)-orbits of all tesselations of ${\mathbb D}$ with doe.  Consider a finitely supported ${\mathbb Z}/2$-{\it marking} $\tau\to\{0,1\}$, and define an equivalence relation on such markings on $\tau$ generated by adding unity modulo two to each of the edges in the frontier of some triangle complementary to $\tau$, i.e., change the values of the marking on all three edges of a fixed complementary triangle.    

Our new universal spin Teichm\"uller space ${\mathcal Tess}^+$ is the collection of all equivalence classes of marked tesselations with doe.   It can be convenient to think of the marking as determined
by finite collections on each edge, where the parity modulo two of the cardinality of the collection determines the
${\mathbb Z}/2$-marking on the edge. 

\begin{figure}[hbt] 
\centering   
\includegraphics[trim =0 430 0 90,width=1.\linewidth]{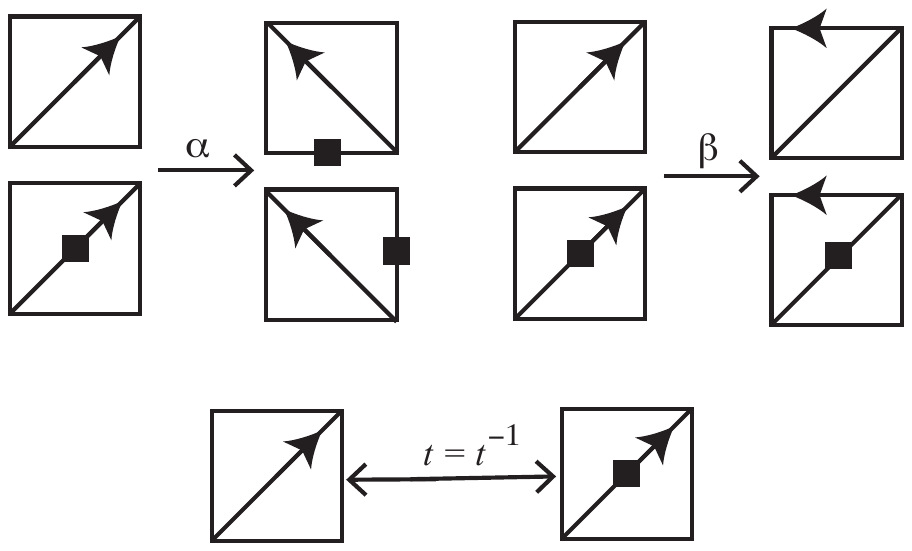}       
\caption{The three generators $\alpha,\beta,t$ of P(SL(2,${\mathbb Z}$)) acting on a marked
tesselation with doe
in the two cases that the doe, which is indicated with an arrow, may or may not have a non-zero marking, as indicated with or without a box icon.
It is easy to confirm that these moves are well-defined on equivalence classes of markings.}
\label{fig:abt}               
\end{figure}

Three combinatorial moves $\alpha,\beta, t$ on marked tesselations with doe are illustrated in Figure \ref{fig:abt}, where adding a mark on an edge corresponds to adding unity modulo two to the value of the marking of the edge.    
The moves on unmarked tesselations with doe
underlying $\alpha$ and $\beta$ generate PPSL(2,${\mathbb Z}$), as we shall prove in the appendix and recall presently, and there is a 
commutative
diagram
$$\begin{array} {c c c c}
{\rm P(SL(2},{\mathbb Z}))\times{\mathcal Tess}^+&\rightarrow&{\mathcal Tess}^+\\
{~~~~~~~~~~~~~~~~~}\downarrow&&\downarrow\\
{\rm PPSL}(2,{\mathbb Z})\times{\mathcal Tess}&\rightarrow&{\mathcal Tess}\\
\end{array}$$
where the vertical maps are the natural forgetful maps.  PPSL(2,${\mathbb Z}$) is conjugate to Thompson's group $T$ \cite{dudes}
by an orientation-preserving  homeomorphism of the circle.

There is a companion \cite{comp} to this paper which explains this isomorphism $T\approx{\rm PPSL}(2,{\mathbb Z})$
and gives the proof that the group 
P(SL(2, ${\mathbb Z}$)) generated by 
$\alpha,\beta, t$  is precisely the group of piecewise-constant maps  $S^1\to$ SL(2,${\mathbb Z}$)
which projectivize to homeomorphisms in PPSL(2,${\mathbb Z}$).   In both these groups, composition is defined by
taking the usual composition on each piece of a common refinement of the two piecewise structures.
The companion \cite{comp} also
surveys the entire construction  of
P(SL(2,${\mathbb Z}$))$\acts{\mathcal Tess}^+$ back to first principles, for punctured surfaces of finite type from \cite{penner} as well
as in the universal case from \cite{unicon}, both of which are covered in \cite{pb}.  
In particular as described there, the new move $\alpha$ including markings arises from the
treatment of spin structures on finite-type punctured surfaces introduced in \cite{N=1}.

This paper derives a finite presentation of the universal spin mapping class group P(SL(2, ${\mathbb Z}$)), for whose 
formulation we require a definition:
Given a word $w=w_1\cdots w_n$ in the generators $\alpha$ and $\beta$, a {\it $t$-insertion} in $w$
is any word $\hat w=w_1 t^{\epsilon_1}w_2 t^{\epsilon _2}\cdots t^{\epsilon_{n-1}} w_n$, for some
choice of $\epsilon _i\in\{0,1\}$.  A further notational point: square brackets in the sequel denote group commutators.

\medskip

\noindent {\bf Main Theorem.}  {\it The group {\rm P(SL(2,${\mathbb Z}$))} generated by $\alpha,\beta,t$ admits a finite presentation with
the following relations}:

\begin{itemize}
\item {\bf Power Laws:} $t^2,\beta^3, \alpha^4, (t\beta)^3, (t\alpha)^4, [t,\alpha^2]~{\rm and}~[t,\alpha t \alpha]$;\\

\item{\bf Pentagon:} $(\beta\alpha)^5~{\rm and}~ (\beta t\alpha t)^5$;\\

\item {\bf Degeneracy:} 

\noindent $
t= [\alpha, \beta t\beta^2]=[\alpha t, \beta t\beta^2]
=[\alpha, t\beta^2t\beta]=[\alpha t, t\beta^2 t\beta];$\\

\item {\bf Insertion:} $[t,\hat w]$ \mbox{\rm ~for any $\mu,\nu\in\{0,1\}$ and any $t$-insertion $\hat w$,}

\noindent {\rm for}~$w\in\{\beta\alpha\beta, ~~
\beta^{2\mu} \alpha^2 ~\beta\alpha\beta ~\alpha^2 \beta^\mu, ~~
\beta^{2\nu}\alpha^2 \beta^{2\mu} \alpha^2~ \beta\alpha\beta~ \alpha^2 \beta^\mu \alpha ^2 \beta ^\nu\}$;\\

\item {\bf First Commutator:} ~\mbox {\rm any word of the form}
$
t^{r_0}\beta\alpha\beta ~\alpha t^{r_4}\alpha ~\beta\alpha\beta$

\noindent $\alpha t^{r_3}\alpha~\beta t^{s_4}\beta t^{s_3}~ \alpha ^3~ t^{s_2} \beta t^{s_1}\beta$ $\alpha t^{r_2} \alpha~\beta t^{t_4}\beta t^{t_3} ~\alpha ^3 ~t^{t_2} \beta t^{t_1}\beta ~\alpha t^{r_1} \alpha,$
{\rm where the exponents satisfy}
$$\begin{aligned}
 \sum_{i=0}^4 r_i&=s_1+s_3+t_1+t_3,\\
r_1+r_2&=s_1+s_2,
\quad r_2+r_3=t_3+t_4,\\
r_3+r_4&=s_3+s_4,
\quad r_1+r_4=t_1+t_2;\\
\end{aligned}$$\\

\item {\bf Second Commutator:} {\rm any word of the form}
$t^{r_0}\beta\alpha\beta ~\alpha t^{r_8}\alpha ~\beta t^{s_5} \beta$

\noindent $\alpha t^{r_7}\alpha~
\beta\alpha\beta~\alpha t^{r_6}~\alpha~\beta t^{s_4}\beta t^{s_3} \alpha ^3
t^{s_2} \beta t^{s_1}\beta 
~\alpha t^{r_4}\alpha\beta t^{t_5}\beta ~\alpha t^{r_3}\alpha~\beta t^{t_4}\beta t^{t_3} \alpha ^3 $

\noindent $t^{t_2} \beta t^{t_1}\beta~\alpha t^{r_2} \alpha\beta~\alpha t^{r_1}\alpha$, 
{\rm where the exponents satisfy} 

$$\begin{aligned}
r_0&=s_1+s_3+s_5+\sum_{i=1}^5 t_i, \\
r_1+r_4&=s_1+s_2,
\quad \hskip -.5exr_2+r_7=t_1+t_2,\\
r_5+r_8&=s_3+s_4,
\quad  r_3+r_6=t_3+t_4,\\
\sum_{i=1}^8 r_i&=s_5+t_1+t_3+t_5.\\
\end{aligned}$$\\
\end{itemize}

 The last two collections of relations arise from $t$-insertions
in the {\it First} $w_1=[\beta\alpha\beta, \alpha^2~\beta\alpha\beta~\alpha^2]$ and
 {\it Second Commutator Relations} $w_2=[\beta\alpha\beta,$ $\alpha^2\beta^2\alpha^2~\beta\alpha\beta~\alpha^2\beta\alpha^2]$ in PPSL(2,${\mathbb Z}$), cf.~Theorem A below.
 There may be redundancies among the Commutator Relations in the Main Theorem.
Note that $\beta$ and $t$ generate a dihedral subgroup $D_6$, 
and $\beta$ and $\alpha^2$ generate a subgroup PSL(2,${\mathbb Z}$) of P(SL(2,${\mathbb Z}$)).
Moreover, $t$ lies in the first derived subgroup according to the Degeneracy Relations, so
P(SL(2,${\mathbb Z}$)) is perfect since $T$ is.

In the First Commutator Relation, the 8-tuple
of $s$- and $t$-variables mimics part of the the root lattice of $E_8$ since $\sum_{i=1}^4 s_i=\sum_{i=1}^4 t_i$.
In the spirit of \cite{FP}, it is natural to wonder if the commutator relations of P(SL(2,${\mathbb Z}$)) 
are organized according to an interesting lattice and 
to ask: What is the automorphism group of P(SL(2,${\mathbb Z}$))?

The basic idea for the proof of the Main Theorem is as follows.
There are decorated  bundles $\widetilde{{\mathcal Tess}}\to {{\mathcal Tess}}$
and $\widetilde{{\mathcal Tess}^+}\to {{\mathcal Tess}^+}$ 
with fibers given by collections of horocycles, one centered at each ideal point
of the tesselation.  As in \cite{unicon}, these respective decorated spaces come equipped
with {\it ideal cell decompositions} ${\mathcal C}$ and ${\mathcal C}^+$, namely, 
decompositions into simplices plus certain of
their faces, where ${\mathcal C}$ and ${\mathcal C}^+$ are
invariant under the respective
actions PPSL(2,${\mathbb Z}$)$\acts\widetilde{{\mathcal Tess}}$
and P(SL(2,${\mathbb Z}$))$\acts \widetilde{{\mathcal Tess}^+}$.

Recall from \cite{unicon} that a codimension-one face in ${\mathcal C}$ corresponds to removing one edge of 
an ideal triangulation.  General position of a path in $\widetilde{\mathcal Tess}$ with respect to ${\mathcal C}$  
therefore shows that
the fundamental path groupoid of $\widetilde{\mathcal Tess}$ is generated by
{\it flips}, that is, the combinatorial moves underlying $\alpha$ in Figure \ref{fig:abt}: namely,
remove an edge from $\tau$ so as to produce a complementary ideal quadrilateral,
and replace the removed edge with the other diagonal of this quadrilateral. 

\begin{figure}[hbt] 
\centering           
\includegraphics[width=0.9\textwidth]{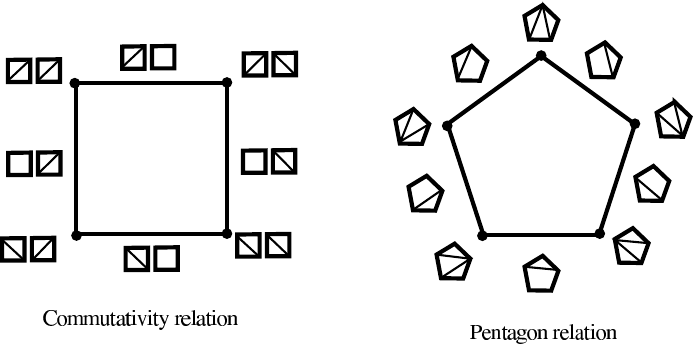}
\caption{Links of codimension-two cells in the ideal cell decomposition of decorated Teichm\"uller spaces.
If one of the edges happens to be the doe, then the pentagon relation has order ten, since the five moves
depicted interchange the two edges, as one can check.} \label{fig:31}               
\end{figure}

General position of a homotopy of paths in $\widetilde{\mathcal Tess}$ likewise shows
that a complete set of relations for the fundamental path groupoid is provided by
the collection of links of codimension-two cells.  These links
correspond to removing two edges, which may either lie in the frontier of a common triangle, or not,
and respectively correspond to the {\it (Classical)  Pentagon} and {\it Commutativity Relations} illustrated in Figure
\ref{fig:31}.
There is a further relation called {\it Idempotence}
arising from a degenerate codimension-two face corresponding to performing a 
flip to produce a new edge upon which you then flip.

It is easy and satisfying to check that ignoring markings, the flip on $\beta(doe)$ is given 
by $\beta\alpha\beta$, and on $\beta ^2(doe)$ by its inverse
$\beta^2\alpha^3\beta^2$.  Instead of flips, the orbifold fundamental group PPSL(2,${\mathbb Z}$)
of the quotient can equally well be regarded as generated by these $\alpha,\beta$, and there is the following
presentation, whose proof is given in the appendix.

\medskip

\noindent {\bf Theorem A.}
{\it  {\rm PPSL(2,${\mathbb Z})$} is generated by the flip $\alpha$ on the doe and the transformation $\beta$ which moves the doe one edge counter-clockwise in the triangle to its left.  A presentation in these generators is given by the following relations:
$\alpha^4$, $\beta^3$, $(\alpha\beta)^5$ and the two commutators $w_1=[\beta\alpha\beta,\alpha^2 \beta\alpha\beta \alpha^2]$ and
$w_2=[\beta\alpha\beta,\alpha^2 \beta^2\alpha^2~\beta\alpha\beta~\alpha^2\beta\alpha^2]$.}

\medskip

\noindent In fact, the current paper provides the first complete proof of this result, which is stated in the context
of \cite{unicon} in \cite{LS} but depends in \cite{LS} upon unpublished computations of Richard Thompson.

Our approach for P(SL(2,${\mathbb Z}$)) is to take each of the defining relations in PPSL(2,${\mathbb Z}$) and consider
all of its $t$-insertions.  As we shall explain presently in an example, we can compute those $t$-insertions which
leave invariant one (and hence each) equivalence class of marking.  These provide a complete but 
highly redundant presentation, which effectively describes the deck group of
$\widetilde{\mathcal Tess}^+/{\rm P(SL(2,}{\mathbb Z}))\to\widetilde{\mathcal Tess}/{\rm PPSL(2,}{\mathbb Z}).$

\begin{figure}[hbt] 
\centering   
\includegraphics[trim =100 410 0 130,width=1.1\linewidth]{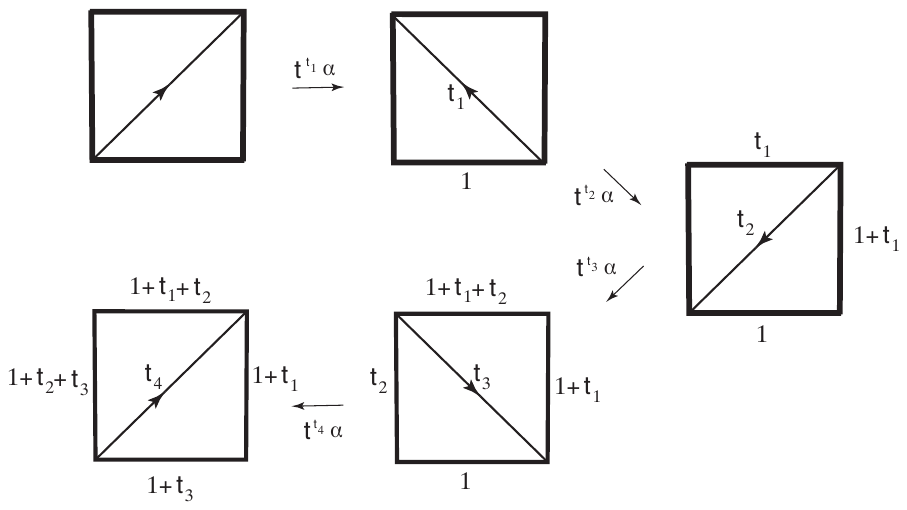}       
\caption{Evolution of marking under the $t$-insertion $\hat w=t^{t_4}\alpha t^{t_3}\alpha t^{t_2}\alpha t^{t_1}\alpha$ in $w=\alpha^4$,
where words act from right to left.}
\label{fig:alphant}               
\end{figure}

\noindent {\bf Example} [The Power Law Relations].  The fourth power of the generator $\alpha\in$~P(SL(2,${\mathbb Z}$)) returns the doe to its starting position and adds a mark on
each edge in the frontier of the quadrilateral near the doe, and likewise for $(t\alpha)^4$. However, this marking is equivalent to the trivial
one with which we started, and so $\alpha^4=1$ in P(SL(2,${\mathbb Z}$)).  More generally, consider 
words $t^{t_4}\alpha t^{t_3}\alpha t^{t_2}\alpha t^{t_1}\alpha$, where
the $t_i\in\{0,1\}$.  We can compute the effects on markings of these words, with letters applied from right to left, as illustrated
in Figure \ref{fig:alphant}.  In order that the resulting marking  is equivalent to the trivial one with which we began, we must have
$t_1=t_3$ and $t_2=t_4$, so we find precisely the Power Law Relations on $\alpha$.
The other Power Law Relations $\beta^3=(t\beta)^3=1$ follow similarly.

\medskip

The overall procedure is thus clear: First, compute the effect on markings of each $t$-insertion in each relation
of PPSL(2,${\mathbb Z}$) to produce a finite presentation of P(SL(2,${\mathbb Z}$)) as in the Example; second, find their redundancies and compute 
a minimal set of relations.  The latter step involves combinatorial group-theoretic calculations, which can be quite involved as in the
Degeneracy Relations, and have yet to be completed for the full presentation in the Main Theorem.

The Insertion Relations were separately discovered in this way, one at a time, and they are in fact special cases of a general lemma: if a word $w$ in $\alpha,\beta, t$ leaves the doe invariant, then $[t,w]=1$.  The general proof is clear, or one can check these several cases directly combinatorially as in the Power Law Example. 

There thus remain four sections, one to analyze each of the
Pentagon, Degeneracy, First and Second Commutator 
Relations, in each case first simplifying as much as possible using the Insertion and Power Law Relations to reduce the number of variables.  The Degeneracy Relation with no $t$-insertions reads simply $(\beta\alpha\beta)(\beta^2\alpha^3\beta^2)$,
which vanishes since $\beta^3=\alpha^4=1$; however with $t$-insertions, there are the four commutators equal to $t$
in the Degeneracy Relations.
The appendix is independent of the body of the paper and constitutes a fifth and final section which derives the finite presentation of
PPSL(2,${\mathbb Z}$) from general position.

Before diving into these detailed computation in  subsequent sections, we collect here several simple 
algebraic facts to be used in the sequel without further comment, throughout which $x,y,t$ are group elements
and $n\in{\mathbb Z}_{>0}$: $(xy)^n=1$ if and only if  $(yx)^n=1$; if $x$ is finite-order, then the following are equivalent:
$[x,y]=1$, $[x^{-1},y]=1$, $ [x,y^{-1}]=1$, $[x^{-1},y^{-1}]=1 $; if $t^2=1$, then $(txty)^n=1$ if and only if $(xtyt)^n=1$, and $(xt)^n=1$ if and only if $(x^{-1}t)^n=1$.

\section{Pentagon Relations}

\begin{figure}[hbt] 
\centering   
\includegraphics[trim =110 550 150 90,width=.8\linewidth]{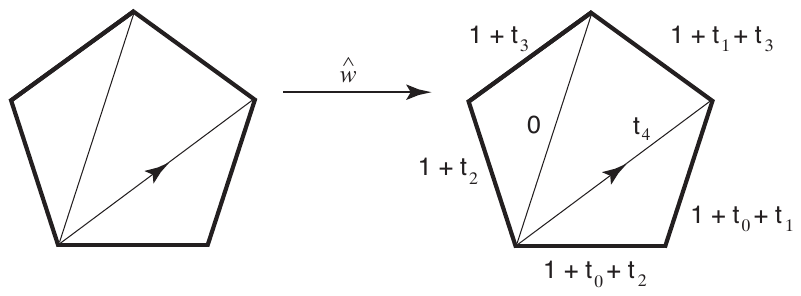}       
\caption{Evolution of marking under the $t$-insertion
$\hat w=\beta\alpha\beta t^{t_4}\alpha\beta t^{t_3}\alpha\beta t^{t_2}\alpha\beta t^{t_1}\alpha t^{t_0}$
in $w=(\beta\alpha)^5$.}
\label{fig:penta}               
\end{figure}

The general insertion of $t$ in $(\beta\alpha)^5$ can be simplified in P(SL(2,${\mathbb Z}$)) 
using the fact from the Insertion Relations that $t$ commutes with any $t$-insertion in $\beta\alpha\beta$.
(In practice, the pentagon $t$-insertion relations provided the vehicle for discovering these Insertion Relations.)
It therefore suffices
to consider only words of the form
$$\hat w=\beta\alpha\beta t^{t_4}\alpha\beta t^{t_3}\alpha\beta t^{t_2}\alpha\beta t^{t_1}\alpha t^{t_0}.$$
Computing as in the Power Law Example, we find the change of marking under $\hat w$ as in Figure \ref{fig:penta},
which is equivalent to the trivial marking if and only if 
$t_1=t_2=t_3$ and $t_4=t_0+t_3$.
The four solutions to this give rise to the relation $(\alpha\beta)^5=1$ as usual,
plus two copies of $(\alpha\beta)^3=t(\alpha\beta t)^3$ and one copy of the familiar $[t,\beta\alpha\beta]=1$.
Meanwhile by the Insertion Relations, $(\alpha\beta)^3=t(\alpha\beta t)^3$ is equivalent to
the relation $(\beta t\alpha t)^5=1$ given in the Main Theorem.

\section{Degeneracy Relations}

\begin{figure}[hbt] 
\centering   
\includegraphics[trim =40 620 150 65,width=.95\linewidth]{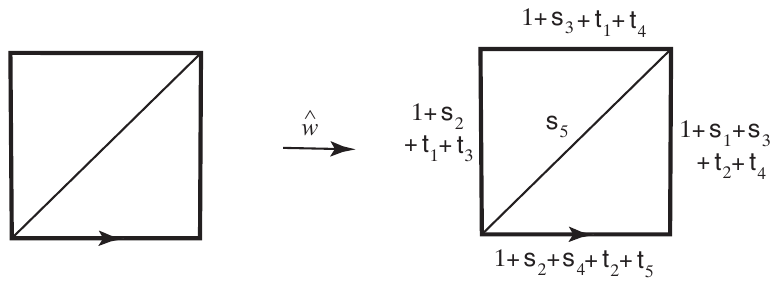}       
\caption{Evolution of marking under the $t$-insertion 
$\hat w= t^{t_5} \beta  t^{s_5} \alpha  t^{t_4} \beta  t^{s_4}\beta  t^{t_3} \beta  t^{s_3}\alpha  t^{t_2}\alpha t^{s_2} \alpha  t^{t_1}\beta t^{s_1} \beta $ in $w=(\beta\alpha\beta)(\beta^2\alpha^3\beta^2)$.}
\label{fig:deg}               
\end{figure}

The action on markings of the $t$-insertion $\hat w$ in $w=(\beta\alpha\beta)(\beta^2\alpha^3\beta^2)$ indicated in Figure \ref{fig:deg}
produces the trivial marking, again as in the Example in the Introduction, if and only if the following linear system holds:
$$\begin{aligned}
t_1+t_2&=s_1+s_5,\\
t_3+t_4&=s_2+s_3,\\
t_3+t_5&=s_1+s_4.\\
\end{aligned}$$
In particular, the complement map that adds one modulo two to each variable preserves solutions,
so to enumerate all $2^{10-3}=128$ solutions, as we shall do here in order to
extract a minimal set of relations, it suffices first to enumerate only those solutions
with at most 5 non-zero variables and then adjoin their complements.  We shall in general let $K$ denote the number of non-zero
variables of a solution.
Define the Boolean predicate
$T(r,s,t)=[(r\wedge s)\vee (s\wedge t)\vee(t\wedge r)]$,
where $\wedge$ is logical AND and $\vee$ is logical OR.
Taking the $\{ 0,1\}$-valued variables as truth values, it is not difficult to combinatorially enumerate the non-zero solutions to the linear system above with $K\leq 5$, as follows:

\medskip

\noindent $\underline{\rm K=2}:~~s_4\wedge t_5,\quad T(s_2,s_3,t_4), \quad T(t_1,t_2,s_5);$

\medskip

\noindent $\underline{\rm K=3}:~~{(i)}~s_1\wedge (s_4\vee t_5) \wedge (t_1\vee t_2\vee s_5),$

\medskip

\noindent\hskip 1.5cm $~~t_3\wedge (s_4\vee t_5) \wedge (s_2\vee s_3\vee t_4)$;

\medskip

\noindent $\underline{\rm K=4}:~~T(s_2,s_3,t_4)\wedge T(t_1,t_2,s_5),$

\medskip

\noindent\hskip 1.5cm $ (ii)~s_4\wedge t_5\wedge 
[T(t_1,t_2,s_5)\vee T(s_2,s_3,t_4)],$;


\medskip

\noindent\hskip 1.5cm $(iii)~s_1\wedge t_3 \wedge (s_2\vee s_3\vee t_4)\wedge (t_1\vee t_2\vee s_5)~$;

\medskip

\noindent $\underline{\rm K=5}:~~s_1\wedge (s_4\vee t_5) \wedge (t_1\wedge t_2\wedge s_5), 
\quad t_3\wedge (s_4\vee t_5) \wedge (s_2\wedge s_3\wedge t_4), $

\medskip

\noindent
\hskip 1.5cm
$s_1\wedge (s_4\vee t_5) \wedge (t_1\vee t_2\vee s_5)\wedge T(s_2,s_3,t_4),$

\medskip

\noindent \hskip 1.5cm 
$t_3\wedge (s_4\vee t_5) \wedge (s_2\vee s_3\vee t_4)\wedge T(t_1,t_2,s_5),$

\medskip

\noindent where we bring cases ($i$)-($iii$) to particular attention.

A tedious group-theoretic computation, which we omit, shows that all but four of
the 128 relations arising from these are tautologies assuming these four new relations together with
the Power Law Relations and the Insertion Relations for $\beta\alpha\beta$.  
Two of the four new relations, $
t\alpha\beta^2t\beta t \alpha^3 \beta t \beta^2$ and $
\alpha t\beta^2t\beta \alpha^3 \beta t \beta^2 t$, 
arise from ($iii$), 
$
\alpha\beta t \beta^2 \alpha^3 t \beta t \beta^2
$
from ($i$), and $
t\beta\alpha t\beta t\beta^2t\alpha^3$
 from ($ii$).
One checks directly  that these equations are identical with the Degeneracy Relations expressed as commutators equal to $t$ as in the Main Theorem.


\section{First Commutator Relations}

\begin{figure}[hbt] 
\centering   
\includegraphics[trim =70 620 110 35,width=1.\linewidth]{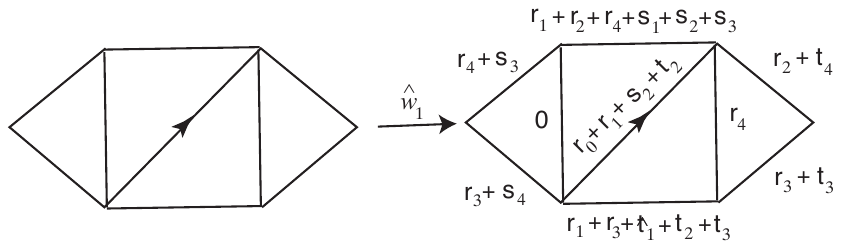}       
\caption{Evolution of marking under the $t$-insertion 
$\hat w_1$ in the First Commutator $w_1$ given the text.}
\label{fig:hex}               
\end{figure}

We use the Insertion Relation for $\beta\alpha\beta$  and the Power Law Relations for $\alpha$ to reduce the dimension of the general $t$-insertion in
the First Commutator $w_1$,
so as to produce

\smallskip

$$\begin{aligned}
\hat w_1= &t^{r_0} \beta\alpha\beta~\alpha  t^{r_4} \alpha~\beta \alpha\beta~\alpha  t^{r_3} \alpha~\\
&\beta  t^{s_4} \beta  t^{s_3} ~\alpha^3  ~t^{s_2}\beta t^{s_1}\beta ~\alpha t^{r_2}\alpha~\\
&\beta  t^{t_4} \beta  t^{t_3} ~\alpha^3  ~t^{t_2}\beta t^{t_1}\beta ~\alpha t^{r_1}\alpha\, .\\
\end{aligned}$$

\smallskip

\noindent The action of $\hat w_1$ on the trivial marking is computed as before with the result illustrated in Figure \ref {fig:hex}, which 
is equivalent to the trivial marking if and only if the five equations in the Main Theorem for the First Commutator are satisfied.

\section{Second Commutator Relations}

We again use the Insertion Relation for $\beta\alpha\beta$ and Power Laws for $\alpha$ to reduce the dimension of the general $t$-insertion in
the Second Commutator $w_2$,
so as to produce

$$\begin{aligned}
\hat w_2=
&t^{r_0}\beta\alpha\beta ~\alpha t^{r_8} \alpha~\beta  t^{s_5} \beta~\alpha  t^{r_7} \alpha~\beta\alpha\beta~\alpha
 t^{r_6} \alpha~\beta~\alpha  t^{r_5} \alpha~\\
 & \beta  t^{s_4} \beta  t^{s_3} ~\alpha^3~  t^{s_2} \beta  t^{s_1} 
 \beta~\alpha  t^{r_4} \alpha~ \beta  \quad t^{t_5} \beta ~\alpha   t^{r_3} \alpha~\\
 & \beta  t^{t_4} \beta  t^{t_3} ~\alpha^3~  t^{t_2} \beta t^{t_1} \beta~\alpha  t^{r_2} \alpha~\beta\quad \alpha  t^{r_1} \alpha\, .\\
  \end{aligned}$$

\begin{figure}[hbt] 
\centering   
\includegraphics[trim =150 600 150 60,width=.85\linewidth]{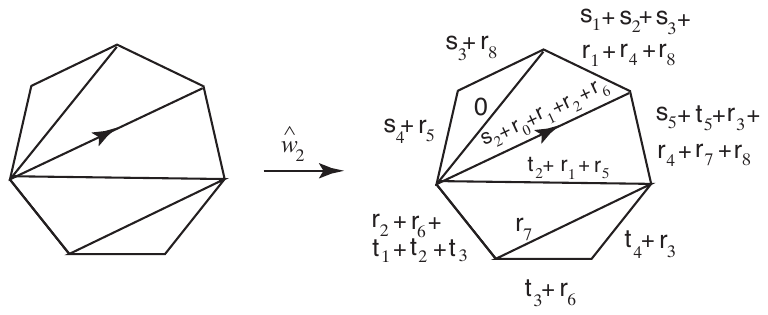}       
\caption{Evolution of marking under the $t$-insertion 
$\hat w_2$ in the Second Commutator $w_2$ given in the text.}
\label{fig:octa}               
\end{figure}

\noindent The action of $\hat w_2$ on the trivial marking is computed as before with the result illustrated in Figure \ref {fig:octa}, which 
is equivalent to the trivial marking if and only if the six equations in the Main Theorem for the Second Commutator are satisfied.

 

\smallskip

\appendix

\section{Presentation of PPSL(2,${\mathbb Z}$)}

This appendix is dedicated to the proof of the following theorem.  The symbols $\alpha,\beta$ here denote
the operations induced from Figure \ref{fig:abt} on tesselations with doe but without marking.

\medskip

\noindent{\bf Theorem A.}~{\it {\rm PPSL(2,${\mathbb Z})$} is generated by the flip $\alpha$ on the doe and the transformation $\beta$ which moves the doe one edge counter-clockwise in the triangle to its left.  A presentation in these generators is given by the following relations:
$\alpha^4$, $\beta^3$, $(\alpha\beta)^5$ and the two commutators $[\beta\alpha\beta,\alpha^2 \beta\alpha\beta \alpha^2]$ and
$[\beta\alpha\beta,\alpha^2 \beta^2\alpha^2~\beta\alpha\beta~\alpha^2\beta\alpha^2]$.}

\begin{proof} According to \cite{unicon},
{\rm PPSL(2,${\mathbb Z})$}  is the orbifold fundamental group of {\rm PPSL(2,${\mathbb Z})$}$\acts\widetilde{\mathcal Tess}$.  This space admits a natural ideal cell
decomposition whose codimension-two skeleton more or less immediately gives this presentation via general position, as discussed in the Introduction.  Namely, the geometry gives three classes of relations:

\medskip

\leftskip .3in

\noindent $\bullet$~the {pentagon relation} $(\alpha\beta)^5=1$ when the two edges cobound a single triangle;

\smallskip

\noindent $\bullet$~{commutativity} of two flips supported on disjoint quadrilaterals when the two edges do not cobound a single triangle;

 \leftskip=0ex
 
 \medskip
 
 \noindent together with the degenerate codimension-two face of consecutively crossing the same codimension-one face:
 
 \leftskip .3in

\medskip

\noindent $\bullet$~{idempotence} of flips, i.e., perform a flip to produce a new edge upon which one subsequently flips.

\medskip\leftskip=0ex

It is evident that the flip $\alpha$ on the doe has order 4, and that $\beta$ has order 3, in the orbifold fundamental group.
Notice that the flip on $\beta(doe)$ is given by $\beta\alpha\beta$ and on $\beta^2(doe)$ is given by
$\beta^2\alpha^3\beta^2$.  
It is well-known that finite words in $\alpha^2\sim \begin{psmallmatrix}0&-1\\1&~0\\\end{psmallmatrix}\in{\rm PSL}(2,{\mathbb Z})$ and $\beta\sim \begin{psmallmatrix}1&-1\\0&~1\\\end{psmallmatrix}\in{\rm PSL}(2,{\mathbb Z})$ act simply transitively on the oriented edges
of the (Farey) tesselation, so the general flip $\phi(g)$ on the edge $g(doe)$, where
$$g=\beta^{\epsilon_n}\alpha^2\beta^{\epsilon_{n-1}}\alpha^2\cdots \alpha^2\beta^{\epsilon _1}
~{\rm with}~\epsilon_j\in\{1,2\}~{\rm for}~n\geq j>1~{\rm and}~\epsilon_1\in\{0,1,2\},$$
can be written 
$$\begin{aligned}
\phi(g)&=g^{-1}\beta^{2\epsilon_n}\alpha^{2\epsilon_n-1}g\\
&=\begin{cases}
\beta^{2\epsilon_1}\alpha^2\cdots\alpha^2\beta^{2\epsilon_{n-1}}\alpha^2~\beta\alpha\beta~
\alpha^2\beta^{\epsilon_{n-1}}\alpha^2\cdots \alpha^2\beta^{\epsilon_1}, {\rm if}~\epsilon_n=1;\\
\beta^{2\epsilon_1}\alpha^2\cdots\alpha^2\beta^{2\epsilon_{n-1}}\alpha^2~\beta^2\alpha^3\beta^2~
\alpha^2\beta^{\epsilon_{n-1}}\alpha^2\cdots \alpha^2\beta^{\epsilon_1}, {\rm if}~\epsilon_n=2,\\
\end{cases}
\end{aligned}$$
applying words from right to left as before.  These are simply the conjugates in PSL(2,${\mathbb Z}$) of the flips
noted above.

In particular, it follows that $\alpha$ and $\beta$ indeed generate PPSL(2,${\mathbb Z}$)
since the group they generate contains the flips on all edges.
Moreover, $[\phi(g)]^{-1}=\phi(\beta^{\epsilon _n}g)$, namely, to invert, change the first (leftmost) exponent
$1\leftrightarrow 2$ of $\beta$ in $g$.
Idempotence of a flip is simply a conjugate of $(\beta\alpha\beta)(\beta^2\alpha^3\beta^2)=1$,
which is not noteworthy for PPSL(2,${\mathbb Z})$ since it follows from $\beta^3=1=\alpha^4$,
but it is of consequence for P(SL(2,${\mathbb Z}$)).
As with flips conjugating by elements of PSL(2,${\mathbb Z}$), the one pentagon relation gives rise to all pentagon relations.

It remains only to prove that the two commutators in the theorem imply the commutativity
relations for flips on any pair of edges which do not bound a common triangle.  This is called a ``remarkable fact'' in \cite{LS}, and here follows its proof by induction:
First note that quite generally two flips
commute if and only if their conjugates commute.  Thus, if flips on two respective edges $e,f$
commute, then flips on any two edges in the same relative positions in the tesselation as $e,f$ also commute.
The first commutator in the statement of the theorem is thus
$[\phi(\beta),\phi(\beta\alpha^2)]=1$, or equivalently 
$$[\alpha^2,\beta\alpha\beta~\alpha^2~\beta\alpha\beta]=1,\leqno{(*)}$$ and this in turn one finds is
also equivalent to $[\phi(\beta),\phi(\beta^2\alpha^2)]=1$; taking inverses one has also
$[\phi(\beta^2),\phi(\beta\alpha^2)]=[\phi(\beta^2),\phi(\beta^2\alpha^2)]=1$.  
Thus, the flip on any edge commutes with the flips on edges {\sl exactly two triangles away},
i.e., a path in general position between the edges meets precisely three edges in the tesselation;
this is equivalent to (*).

The flip $\alpha$ on the doe does not commute with other flips because it engenders a change of 
notation, cf. Figure \ref{fig:alpha}.  Rather, the flip $\phi(g)$ on $g(doe)$ evidently satisfies the following functional equation
$$\hskip -8em{\rm FE}(g):\quad\quad\quad\hskip 5em\alpha~\phi (g)=\phi(T(g))~\alpha,$$
where in the notation for $g$ above 
$$T(g)=\begin{cases}
g\beta^{\epsilon_1}\alpha^{2\epsilon _1},&{\rm if}~\epsilon_1\neq 0;\\
g\alpha^2\beta^{\epsilon_2}\alpha^{2\epsilon _2},&{\rm if}~\epsilon_1=0,\\
\end{cases}$$
or in other words the transformation $T$ acts on suffixes of $g$ by
$$\beta\mapsto\beta^2\alpha^2\mapsto\beta\alpha^2\mapsto\beta^2\mapsto\beta.$$
In particular, $T$ is a square root of right multiplication by $\alpha^2$.

\begin{figure}[hbt] 
\centering   
\includegraphics[trim =0 490 0 120,width=1.\linewidth]{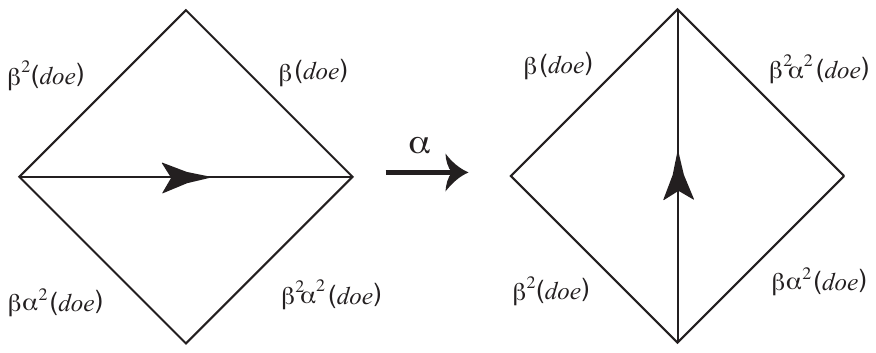}       
\caption{Change of notation under $\alpha$ near the doe.}
\label{fig:alpha}               
\end{figure}

The flip on the edge underlying the doe commutes with the flip $\phi(g)$
if and only if $g$ satisfies the functional equation FE($g$).   Thus, the first
commutator in the statement of the theorem is equivalent to FE($g$) holding for all $g=\beta^{\epsilon_2}\alpha^2 \beta^{\epsilon_1}$
with $\epsilon_1\neq 0$ and $n=2$ in our normal form. 
Moreover if $\epsilon_1(g)\neq 0$, then FE($g$) trivially implies FE($g\alpha^2$), so it suffices to assume
henceforth that $\epsilon_1\neq 0$ and proceed by induction on $n$, as we next undertake.
One confirms that the case $n=3$ of FE is likewise equivalent to the identity
$$[v,\beta v \beta^2]=1, ~{\rm where}~v=\alpha^2 ~\beta\alpha\beta ~\alpha^2,\leqno{(\dag)}$$
which is evidently equivalent to $[\beta^{\epsilon}v^{\pm 1}\beta^{2\epsilon}, \beta^{\delta}v^{\pm 1}\beta^{2\delta}]=1$, for any
$\epsilon, \delta\in\{ 1,2\}$, and finds that 
($\dag$) is in turn equivalent to the vanishing of the second commutator in the statement of the theorem.

For the inductive step, suppose that $g=\beta^{\epsilon_n}\alpha^2\cdots \alpha^2 \beta^{\epsilon_1}$ with $n>3$ in our normal form with $\epsilon_1=\epsilon_1(g)\neq 0$, 
and set $h=\beta^{\epsilon_{n}}\alpha^2\cdots \alpha^2 \beta^{\epsilon_2}$.  FE($g$) is given by
$$\begin{aligned}
\alpha ~ \phi(h\alpha^2\beta^{\epsilon _1})&= \phi(h \alpha^2\beta^{\epsilon _1}~ \beta^{\epsilon _1} \alpha ^{2\epsilon _1})~\alpha\\
&=\alpha^{2\epsilon _1}~\phi(h\alpha^2 \beta^{2\epsilon_1})~\alpha^{2\epsilon_1+1},\\
\end{aligned}$$
so for $\epsilon_1(g)\neq 0$, FE($g$) reads
$$\alpha^{2\epsilon _1+1}~\phi(g) ~\alpha^{2\epsilon _1-1}= \beta^{2\epsilon_1}~\phi(g)~\beta^{\epsilon _1},$$
which simply means: 
 if $\epsilon_1(g)=2$, then you can pull $\alpha$ to the right across $\phi(g)$
at the expense of changing this to $\epsilon_1 =1$, or equivalently to the left  again changing the terminal exponent if $\epsilon_1(g)=1$.
 In particular for any $N\geq 1$, FE($g$) holds for $g=(\beta\alpha^2)^N\beta$ and likewise for $\beta g$.  
 
Suppose first that in fact $\epsilon_1(g)=1$. If $g$ differs from the forms above that automatically satisfy FE, then
there is some index $1<m<n$ with $\epsilon _m =2$, so that $g=h\alpha^2 k$ where $h$ ends with $\beta^2$ and $k$ of course
ends with $\beta$ since we assume that $g$ does, and both $h$ and $k$ satisfy FE by the strong inductive hypothesis.  Using the simple
description of the FE as right/left commutativity laws for $\alpha$ across flips in the previous paragraph, one finds that FE($g$) is equivalent to $[\phi(h'),\phi(k')]=1$,
where $h',k'$ respectively arise from $h,k$ by altering their terminal $\beta$-exponents.
 Since the relative positions of the flipped edges
is decreased by unity, FE holds in general by induction.  The analogous argument holds for $\epsilon_1(g)=2$ using the automatic solutions
${\rm FE}(\beta^2(\alpha^2\beta^2)^N \alpha^2)$, for $N\geq 1$. \end{proof}

\medskip

Notice that according to the proof there are exactly two conjugacy classes in ${\rm PPSL}(2,{\mathbb Z})$ of flips on edges other than the doe, namely,
the conjugacy classes of $\beta\alpha\beta$ and $\beta^2\alpha^3\beta^2$.  Since these are inverses,  the collection of flips on edges other than the doe abelianizes to a cyclic group.  Inspection of Figure \ref{fig:31} shows that the five flips of the pentagon relation  are comprised of
two from one class and three from the other.  It follows that the pentagon relations
alone imply that flips on edges other than the doe lie in the first derived subgroup.

\vfill\eject

\end{document}